\documentstyle[avv3,12pt]{article}

\newcounter{minutes}
\divide\time by 60
\newcounter{hours}
\multiply\time by 60
\addtocounter{minutes}{-\time}
\begin{document}
% Macros for AVV book

%\def\R{$\rm{I}\!\rm{R}$} Ei kelpaa
%\def\N{$\rm{I}\!\rm{N}$}  -''-

\def\wrsty{\bf\large}
\def\sq{\ \ \square }
\def\one{{$(1)\ $}}

\def\R{\,\brm {R}}
\def\N{\brm {N}}
\def\C{\brm {C}}
\def\Z{\brm {Z}}

\def\B{\,\brm {B}} 

\def\Bcal{\cal {B}}
\def\Bn{\B ^{n}}
\def\Bo{\overline{\B}}
\def\Bon{\overline{\B} ^n}
\def\sfm{\msf {M}}

\def\H{\brm {H}}

\def\Hn{\H ^{n}}

\def\rdt{{r'}^{\,2}}
\def\xdt{{x'}^{\,2}}
\def\rdth{{r'}^{\,3}}
\def\rdf{{r'}^{\,4}}
\def\sdt{{s'}^{\,2}}

\def\prkinf{\prod_{k=1}^\infty}
\def\prninf{\prod_{n=1}^\infty}
\def\prkn{\prod_{k=1}^n}
\def\prkp{\prod_{k=1}^p}
\def\prkq{\prod_{k=1}^q}
\def\prkNn{\prod_{k=N}^n}

\def\AG{\mathop{\rm{AG}}\nolimits}
\def\G{\mathop{\rm{G}}\nolimits}
\def\A{\mathop{\rm{A}}\nolimits}

% \L changed to \La  5.12.94 because \L is needed for polish L
%\def\La{\rm {L}\,}
\def\La{\mathop{\rm{L}}\nolimits}

% muutin \def\th{\rm {th}\,} -> \def\th{\mathop{\rm{th}}\nolimits}
% jne., P. S. 2.12.94
\def\th{\mathop{\rm{th}}\nolimits}
\def\sh{\mathop{\rm{sh}}\nolimits}
\def\sn{\mathop{\rm{sn}}\nolimits}
\def\ns{\mathop{\rm{ns}}\nolimits}
\def\nc{\mathop{\rm{nc}}\nolimits}
\def\nd{\mathop{\rm{nd}}\nolimits}
\def\sd{\mathop{\rm{sd}}\nolimits}
\def\cd{\mathop{\rm{cd}}\nolimits}
\def\cs{\mathop{\rm{cs}}\nolimits}
\def\ds{\mathop{\rm{ds}}\nolimits}
\def\dc{\mathop{\rm{dc}}\nolimits}
\def\ch{\mathop{\rm{ch}}\nolimits}
\def\cn{\mathop{\rm{cn}}\nolimits}
\def\tn{\mathop{\rm{tn}}\nolimits}
\def\dn{\mathop{\rm{dn}}\nolimits}
\def\csch{\mathop{\rm{csch}}\nolimits}
\def\sech{\mathop{\rm{sech}}\nolimits}
\def\coth{\mathop{\rm{coth}}\nolimits}
\def\arth{\mathop{\rm{arth}}\nolimits}
\def\arsh{\mathop{\rm{arsh}}\nolimits}
\def\arch{\mathop{\rm{arch}}\nolimits}
\def\arcsech{\mathop{\rm{arcsec}}\nolimits}
\def\Arc{\mathop{\rm{Arc}}\nolimits}
\def\Arg{\mathop{\rm{Arg}}\nolimits}

\def\Log{\mathop{\rm{Log}}\nolimits}

\def\a{\alpha}

\def\eqb{\begin{equation}}
\def\eqe{\end{equation}}
\def\eb{\begin{eqnarray}}
\def\ee{\end{eqnarray}}
\def\ebnn{\begin{eqnarray*}}
\def\eenn{\end{eqnarray*}}
\def\db{\begin{displaystyle}}
\def\de{\end{displaystyle}}
\def\tb{\begin{textstyle}}
\def\te{\end{textstyle}}
\def\exb{\begin{ex}}
\def\exe{\end{ex}}

%%%%%%%%% KAUNIS K \K  %%%%%%%%%%%%%%%%%%%%%%%%%%%%%%%%%%%%%%%%%%%%%
\font\fFt=eusm10 scaled 1200
\font\fFa=eusm7 scaled 1200
\font\fFp=eusm5 scaled 1200
\def\K{\mathchoice
%displaystyle
{\hbox{\,\fFt K}}
%textstyle
{\hbox{\,\fFt K}}
%scriptstyle
{\hbox{\,\fFa K}}
%scriptscriptstyle
{\hbox{\,\fFp K}}}

%%%%%%%%% KAUNIS E \E  %%%%%%%%%%%%%%%%%%%%%%%%%%%%%%%%%%%%%%%%%%%%%

\def\E{\mathchoice
%displaystyle
{\hbox{\,\fFt E}}
%textstyle
{\hbox{\,\fFt E}}
%scriptstyle
{\hbox{\,\fFa E}}
%scriptscriptstyle
{\hbox{\,\fFp E}}}
%%%%%%%%%%%%%%%%%%%%%%%%%%%%%%%%%%%%%%%%%%%%%%%%%%%%%%%%%%%%%%%%%%%%%

%%%%%%%%% KAUNIS R \Rb  %%%%%%%%%%%%%%%%%%%%%%%%%%%%%%%%%%%%%%%%%%%%%

\def\Rb{\mathchoice
%displaystyle
{\hbox{\fFt R}}
%textstyle
{\hbox{\fFt R}}
%scriptstyle
{\hbox{\fFa R}}
%scriptscriptstyle
{\hbox{\fFp R}}}
%%%%%%%%%%%%%%%%%%%%%%%%%%%%%%%%%%%%%%%%%%%%%%%%%%%%%%%%%%%%%%%%%%%%%

%%%%%%%%% KAUNIS F \Fb  %%%%%%%%%%%%%%%%%%%%%%%%%%%%%%%%%%%%%%%%%%%%%

\def\Fb{\mathchoice
%displaystyle
{\hbox{\fFt F}}
%textstyle
{\hbox{\fFt F}}
%scriptstyle
{\hbox{\fFa F}}
%scriptscriptstyle
{\hbox{\fFp F}}}
%%%%%%%%%%%%%%%%%%%%%%%%%%%%%%%%%%%%%%%%%%%%%%%%%%%%%%%%%%%%%%%%%%%%%

\def\Ma{{\fFt A}}
\def\Mb{{\fFt B}}
\def\Mc{{\fFt C}}
\def\Md{{\fFt D}}
\def\Me{{\fFt E}}
\def\Mf{{\fFt F}}
\def\Mg{{\fFt G}}
\def\Mh{{\fFt H}}
\def\Mi{{\fFt I}}
\def\Mj{{\fFt J}}
\def\Mk{{\fFt K}}
\def\Ml{{\fFt L}}
\def\Mm{{\fFt M}}
\def\Mn{{\fFt N}}
\def\Mo{{\fFt O}}
\def\Mp{{\fFt P}}
\def\Mq{{\fFt Q}}
\def\Mr{{\fFt R}}
\def\Ms{{\fFt S}}
\def\Mt{{\fFt T}}
\def\Mu{{\fFt U}}
\def\Mv{{\fFt V}}
\def\Mw{{\fFt W}}
\def\Mx{{\fFt X}}
\def\My{{\fFt Y}}
\def\Mz{{\fFt Z}}

\def\cc{\clearpage\setcounter{equation}{0}
\setcounter{figure}{0}\setcounter{table}{0}}

%****************************************

\def\q#1{_{#1}^{}}

\def\cc{\setcounter{equation}{0}}
\def\A{{\cal A}}
\def\B{{\cal B}}  % m.tex:in \B ei ole kaytossa enaan
\def\F{{\cal F}}
\def\J{{\cal J}}
\def\M{{\cal M}}
\def\N{{\cal N}}
\def\O{{\cal O}}
\def\U{{\cal U}}
\def\WT{{\cal W}{\cal T}}
\def\X{{\cal X}}
\def\Y{{\cal Y}}

\def\1{{\hbox{\rm 1}\hskip-0.38em\hbox{\rm 1}}}

%\def\C{\Bbb C} oli 10pt tekstia, tehdaan uusi:
%\font\bbb=msym10 scaled 1200
\def\C{\mbox{\bf C}}
\def\isowedge{{\textstyle\bigwedge\limits}}
\def\Re{{\,\rm Re}}
\def\Im{{\,\rm Im}}
\def\osc{{\rm osc}}
\def\inter{{\rm int}}
\def\supp{{\rm supp}\,}
\def\loc{{\rm loc}}
\def\dist{{\rm dist}\,}
\def\diver{{\rm div}\,}
\def\epsilon{\varepsilon}

\newcommand{\vse}{\vspace{.2in}}
\newcommand{\np}{\newpage}

\np
\pagenumbering{arabic}
\pagestyle{myheadings}
\markboth{\today}{\today}

\begin{center}
{\Large\bf Ahlfors theorems for differential forms}
\end{center}
\medskip

\begin{center}
{\large\bf O.~Martio, V.M.~Miklyukov,  and M.~Vuorinen}
\end{center}
\bigskip
\centerline{\tt File: 
~\jobname .tex, printed: \number\year-\number\month-\number\day,
        \thehours.\ifnum\theminutes<10{0}\fi\theminutes}
\medskip

{\centerline{\bf Abstract}

Some counterparts of theorems of
Phragm\'en-Lindel\"of and of Ahlfors are proved
for differential forms of ${\cal WT}$--classes.}

\bigskip

\section{$\cal WT$-forms}{}
%Our notation is as in \cite{FMMVW}. We recall some results on differential 
%forms. For connections with qr-mappings see {\cite{FMMVW}}. 
This paper is continuation of the earlier work \cite{FMMVW}, where the
main topic was to examine the connection between quasiregular (qr) mappings 
and so called ${\cal WT}$--classes of differential forms. 
We first recall some basic notation and 
terminology from {\cite{FMMVW}}. 

Let $\M$ be a Riemannian manifold of class $C^3$, $\dim \M=n$, with or
without boundary, and let
\eqb
w \in L_{\loc}^p(\M),\ \deg w=k,\ 0\le k\le n,\
p>1,\label{eq2.2}
\eqe
be a weakly closed differential form on $\M$, i.e. for each form
$$
\varphi\in W^1_{q,{\rm loc}}(\M),\quad {\rm deg}\,\varphi=k+1,
\quad {1\over p}+{1\over q}=1,
$$ 
with a compact ${\rm supp}\,\varphi$ in $\M$ and such that 
${\rm supp}\,\varphi\cap \partial \M=\emptyset$, we have
$$
\int\limits_{\M}\langle w,\delta\varphi\rangle\,dv=0\,.
$$ 
Here $\delta\varphi=(-1)^k\star^{-1} d\star\varphi$, $k={\rm deg}\,\varphi$
and $\star\alpha$ is the orthogonal complement of a differential form $\alpha$ 
on a Riemannian manifold $\M$.

A weakly closed form $w$ of the kind (\ref{eq2.2}) is said to be of the class
$\WT_1$ on $\M$ if there exists a  weakly closed differential form
\eqb
\theta\in L_{\loc}^q(\M),\ \deg \theta=n-k,\
{1\over p}+{1\over q}=1,\label{eq2.4}
\eqe
such that almost everywhere on $\M$ we have
\eqb
{\nu}_0\ |\theta|^q\le \ \langle w,*\theta\rangle\label{eq2.5}
\eqe
for some constant ${\nu}_0$.
\smallskip

The differential form (\ref{eq2.2}) is said to be of the class $\WT_2$
on $\M$ if there exists a differential form
(\ref{eq2.4}) such that almost everywhere on $\M$ 
\eqb
\nu_1\,|w|^p\le \ \langle w,*\theta\rangle\label{eq2.7}
\eqe
and
\eqb
|\theta|\le \nu_2\,|w|^{p-1}\label{eq2.8}
\eqe
for some constants ${\nu}_1,{\nu}_2>0$.
\smallskip

\begin{thm}{}\label{2.16}
$\WT_2\subset\WT_1.$
\end{thm}
\medskip

For a proof see \cite{FMMVW}.
\smallskip

The following partial integration formula for differential forms is 
useful \cite{FMMVW}.
\smallskip

\begin{lem}{}\label{1.33a}
Let $\alpha\in W^1_{p,\loc}(\M)$ and $\beta\in W^1_q(\M)$ be differential 
forms, $\deg\alpha+\deg\beta=n-1$, 
${1/p}+{1/q}=1, \quad 1\le p,\,q\le \infty$, and let 
$\beta$ have a compact support ${\rm supp}\,\beta \subset {\M}$. Then
\eqb
\label{eq1.34a}
\int\limits_\M d\alpha\wedge\beta=(-1)^{\deg\alpha+1}\int\limits_\M
\alpha\wedge d\beta.
\eqe
In particular, the form $\alpha$ is weakly closed if and only if 
$d\alpha=0$ a.e.\ on $\M$.
\end{lem}
\smallskip

Let $\A$ and $\B$ be Riemannian manifolds of dimensions $\dim
\A=k$, $\dim\B=n-k$, $1\le k<n$, and with scalar products
$\langle\, ,\rangle_A$, $\langle\, ,\rangle_B$, respectively. The Cartesian
product $\N=\A\times\B$ has the natural structure of a Riemannian
manifold with the scalar product
$$\langle\, ,\rangle=\langle\, ,\rangle_{\A}+\langle\, ,\rangle_{\B}.$$
We denote by $\pi:\A\times\B\to\A$ and $\eta:\A\times\B\to\B$ the
natural projections of the manifold $\N$ onto submanifolds.

If $w_{\A}$ and $w_{\B}$ are volume forms on $\A$ and $\B$, respectively,
then the differential form $w_{\N}=\pi^*w_{\A}\wedge\eta^*w_{\B}$ 
is a volume form on $\N$.
\smallskip

Let $y_1,\ldots,y_k$ be an orthonormal system of coordinates
in ${\bf R}^k$, $1\le k\le n$. Let $\A$ be a domain in 
${\bf R}^k$ and let $\B$ be
an $(n-k)$-dimensional Riemannian manifold. We consider the manifold
$\N=\A\times\B$.

\bigskip

%%%SECTION
%%%SECTION
%%%SECTION
\cc
\section{Boundary sets}{}\label{sec3}
Below we introduce the notions of parabolic and hyperbolic type of
boundary sets on noncompact Riemannian manifolds and study exhaustion
functions of such sets. We also present some illuminating examples.
\bigskip

%\begin{subsection}{Boundary sets.}\label{3.1}
%\end{subsection}
 Let $\M$ be an $n$-dimensional noncompact
Riemannian manifold without boundary. Boundary sets on $\M$
are analogies to prime ends
due to Caratheodory (cf.\ e.g.\ \cite{Su}).

Let $\{\U_k\}$, $k=1,2,\ldots$ be a collection of open sets
$\U_k\subset\M$ with the following properties:

(i) for all $k=1,2,\ldots\quad$ $\overline{\U}_{k+1} \subset\U_{k}$,

(ii) $\bigcap\limits_{k=1}^{\infty}  \overline{\U}_k=\emptyset$.

\smallskip

A sequence with these properties will be called a chain on the
manifold $\M$.
\smallskip

Let $\{\U'_k\}$, $\{\U''_k\}$ be two chains of
open sets on $M$. We shall say that the chain $\U'_k$ is {\it contained}
in the chain $\{\U''_k\}$, if for each $m\ge 1$ there exists a number
$k(m)$ such that for all $k>k(m)$ we have $\U'_k\subset\U''_m$. Two
chains, each of which is contained in the other one, are called {\it
equivalent}. Each equivalence class $\xi$ of chains is called a {\it
boundary set} of the manifold $\M$. To define $\xi$ it is enough to
determine at least one representative in the equivalence class. If the
boundary set $\xi$ is defined by the chain $\{\U_k\}$, then we shall
write $\xi\asymp\{\U_k\}$.
\smallskip

A sequence of points $m_k\in\M$ converges to $\xi$ if for
some (and, therefore, all) chain $\{\U_k\}\in\xi$ the following
condition is satisfied: for every $k=1,2,\ldots$ there exists an integer
$n(k)$ such that $m_n\in\U_k$ for all $n>n(k).$  A sequence $(m_n)$
lies off a boundary set $\xi\asymp\{\U_k\}$, if for every $k=1,2,\ldots$
there exists a number $n(k)$ such that for all $n>n(k)\ m_n\notin\U_k$.

\smallskip

A boundary set
$\xi\asymp\{\U_k\}$ is called a set of ends of the manifold $\M$ if each of
$\{\U_k\}$ has a compact boundary $\partial\U_k$. If in addition each of 
the sets $U_k$ is connected, then $\xi\asymp\{U_k\}$ is called 
{\it an end of the manifold} $\M$.
\bigskip 

\begin{nonsec}{Types of boundary sets }\label{3.5}
\end{nonsec}
 Let $D$ be an open set on $\M$ and let
$A,B\subset D$ be closed subsets in $D$ such that 
$\overline A\cap\overline B=\emptyset$. 
Each triple $(A,B;D)$ is called a condenser on $\M$.

We fix $p\ge 1$. The $p$-capacity of the condenser $(A,B;D)$ is defined
by
\eqb
\hbox{\rm cap}_p(A,B;D)=\inf\int\limits_D 
|\nabla\varphi|^p*\1_{\M},\label{eq3.6}
\eqe
where the infimum is taken over the set of all continuous functions of
class $W_{p,\loc}^1(D)$ such that $\varphi(m)|_A=0$,
$\varphi(m)|_B=1$. It is easy to see that for a pair $(A,B;D)$ and
$(A_1,B_1;D)$ with $A_1\subset A$, $B_1\subset B$ we have
$$
\hbox{\rm cap}_p(A_1,B_1;D)\le\hbox{\rm cap}_p(A,B;D).
$$

Let $\overline{A}$, $\overline{B}$ be compact in $D$. A standard 
approximation method shows that 
$\hbox{\rm cap}_p(A,B;D)$ does not change if one restricts the class
of functions in the variational problem (\ref{eq3.6}) to Lipschitz
functions equal to $0$ and $1$ in the sets $A$ and $B,$ respectively.

\smallskip

Let $\{\U_k\}$ be an arbitrary chain on a
manifold $\M$. We fix a subdomain $H\subset\subset\M$. If $k$ is
sufficiently large, the intersection $\overline
H\cap\overline{\U}_k=\emptyset$ and we consider the
condenser $(\overline H,\overline{\U_k};\M)$. Then it is clear that for
$k=1,2,\ldots$
$$\hbox{\rm cap}_p(\overline H,\overline{\U_k};\M)\ge\hbox{\rm
cap}_p(\overline H,\overline{\U_{k+1}};\M).$$
We shall say that the chain $\{\U_k\}$ on $\M$ has $p$-{\it capacity
zero}, if for every subdomain $H\subset\subset\M$ we have
\eqb
\lim\limits_{k\to\infty}\hbox{\rm cap}_p (\overline
H,\overline{\U_k};\M)=0.\label{eq3.8}
\eqe

We shall say that a boundary set $\xi$ is of $p$-{\it parabolic type} if
every chain $\{\U_k\}\asymp\xi$ is of $p$-capacity zero. A boundary set
$\xi$ is of $\alpha$-{\it hyperbolic type} if at least one of the chains
$\{\U_k\}\in\xi$ is not of $p$-parabolic type.
\smallskip

%\begin{subsec}{}\label{3.9}
%\end{subsec}
Let 
$$
\{\U_k\}_{k=1}^{\infty},\quad \overline{\U}_k\subset \U_{k+1},
\quad \bigcup\nolimits_{k=1}^{\infty}\U_k=\M
$$ 
be an arbitrary exhaustion of
the manifold $\M$ by subdomains $\{\U_k\}$. The manifold $\M$ is of
$p$-parabolic or $p$-hyperbolic type depending on the
$p$-parabolicity or $p$-hyperbolicity of the boundary set
$\{\M\setminus\overline{\U}_k\}$.

It is well--known, see \cite{KE}, that a noncompact Riemannian manifold
$\M$ without boundary is of $p$-parabolic type if and only if
every solution of the inequality
$$
\diver_{\M}(|\nabla u|^{p-2}\nabla u)\ge 0
$$
which is bounded from above is a constant.

The classical parabolicity and hyperbolicity coincides with $2$-parabolicity 
and $2$-hyperbolicity,
respectively. Therefore whenever we refer to
parabolic or hyperbolic type (of a manifold or a boundary
set) we mean $2$-parabolicity or $2$-hyperbolicity.
\smallskip

\begin{exmp}\label{3.10}
The space ${\bf R}^n$ is of $p$-parabolic type for $p\ge n$ and $p$-hyperbolic
type for $p<n$.
\end{exmp}
\smallskip

We next present a proposition that provides a convenient method of
verifying the $p$-parabolicity and $p$-hyperbolicity of boundary sets.

\smallskip

\begin{lem}{$\;$\cite{MIK3}}\label{3.11}
 Let $\xi$ be a boundary set on $\M$. If for a chain
$\{\U_k\asymp\xi\}$ and for a nonempty open set $H_0\subset\subset\M$ the
condition (\ref{eq3.8}) holds, then the boundary set $\xi$ is of $p$-parabolic
type.
\end{lem}
\smallskip
\begin{nonsec}{$A$-solutions }
\end{nonsec}
Let $\M$ be a Riemannian manifold and let
$$
A: T(\M) \to T(\M)
$$ be a
mapping defined a.e. on the tangent bundle $ T(\M) .$
Suppose that for a.e. $m \in \M$ the mapping $A$ is continuous on the fiber
$T_m , $ i.e. for a.e. $m \in \M$ the function
$ A(m, \cdot ): \xi \in T_m  \to T_m $ is defined and continuous; the
mapping $m \to A_m (X)$ is measurable for all measurable vector fields
$X$ (see \cite{HKM}).

Suppose that for a.e. $m \in \M$ and for all 
$\xi\in T_m$ the inequalities  
\eqb
\nu_1\,|\xi|^p\le \langle\xi, A(m,\xi)\rangle,\label{eq2.23}
\eqe
and
\eqb
|A(m,\xi)|\le \nu_2\,
|\xi|^{p-1}\label{eq2.24}
\eqe
hold with $p>1$ and for some constants $\nu_1,\nu_2>0$. It is clear that 
we have $\nu_1\le\nu_2$.

We consider the equation  
\eqb
\diver\,A(m,\nabla f)=0.\label{eq2.19}
\eqe
Solutions to (\ref{eq2.19}) are understood in the weak sense, 
that is, $A$-solutions are $W_{p,loc}^1$-functions satisfying the integral
identity
\eqb
\label{eq2.20}
\int\limits_{\M}\langle\nabla\theta,A(m,\nabla f)\rangle *\1_{\M}=0
\eqe
for all $\theta\in W_p^1(\M)$ with a compact support ${\rm theta}\subset{\M}$.

A function $f$ in $W_{p,loc}^1 (\M)$ is a {\it $A$-subsolution} of 
(\ref{eq2.19}) in $\M$ if
\eqb
\diver\,A(m,\nabla f)\ge 0\label{eq2.25}
\eqe
weakly in $\M$, i.e.
\eqb
\label{eq2.26}
\int\limits_{\M}\langle\nabla\theta,A(m,\nabla f)\rangle *\1_{\M}\le 0
\eqe
whenever $\theta\in W^1_p (\M)$, is nonnegative with 
compact support in  $\M$.
\medskip

A basic example of such an equation is the $p$-Laplace equation
\eqb
\label{ieq3.21}
\diver(|\nabla f|^{p-2}\nabla f)=0.
\eqe
\bigskip

%%%SECTION
%%%SECTION
%%%SECTION

\cc
\section{Exhaustion functions}{}
Below we introduce exhaustion and special exhaustion functions on 
Riemannian manifolds and give illustrating examples.
\bigskip

\begin{nonsec}{Exhaustion functions of boundary sets }\label{3.14}
\end{nonsec}
Let $h:\M\to(0,h_0)$, $0<h_0\le\infty$, be a locally Lipschitz function. For
arbitrary $t\in(0,h_0)$ we denote by
$$
B_h(t)=\{m\in\M:h(m)<t\},\quad
\Sigma_h(t)=\{m\in\M:h(m)=t\}
$$
the $h$-balls and $h$-spheres, respectively.
\smallskip

Let $h:\M\to{\bf R}$ be a locally Lipschitz function such that: 
there exists a compact $K\subset\M$ such that $|\nabla h(x)|>0$ 
for a.e. $m \in \M \setminus K$.
We say that the function $h$ is an exhaustion function for a boundary
set $\Xi$ of $\M$ if for an arbitrary sequence of points $m_k\in\M$,
 $k=1,2,\ldots$ the function $h(m_k)\to h_0$ if and only if $m_k\to\xi$.

It is easy to see that this requirement is
satisfied if and only if for an arbitrary increasing sequence
$t_1<t_2<\ldots<h_0$ the sequence of the open sets
$V_k=\{m\in\M:h(m)>t_k\}$ is a chain, defining a boundary set $\xi$.
Thus the function $h$ exhausts the boundary set $\xi$ in the traditional sense
of the word. 
\smallskip

The function $h:\M\to(0,h_0)$ is called the exhaustion function of the
manifold $\M$ if the following two conditions are satisfied

(i) for all $t\in(0, h_0)$ the $h$--ball $\overline{B_h(t)}$ is compact;

(ii) for every sequence $t_1<t_2<\ldots < h_0$ with
$\lim\nolimits_{k\to\infty}t_k=h_0$, the sequence of $h$-balls
$\{B_h(t_k)\}$ generates an exhaustion of $\M$, i.e. 
$$
B_h(t_1)\subset B_h(t_2)\subset\ldots\subset B_h(t_k)\subset\ldots 
\quad \mbox{and} \quad
\cup_k B_h(t_k)=\M.
$$
\smallskip

\begin{exmp}\label{3.17}
Let $\M$ be a Riemannian manifold. We set
$h(m)=\mbox{dist}(m,m_0)$ where $m_0\in\M$ is a fixed point. Because
$|\nabla h(m)|=1$ almost everywhere on $\M$, the function $h$
defines an exhaustion function of the manifold $\M$.
 \end{exmp}

\bigskip

\begin{nonsec}{Special exhaustion functions }\label{3.20}
\end{nonsec}
Let $\M$ be a noncompact Riemannian manifold with the boundary $\partial \M$
(possibly empty).
Let $A$ satisfy (\ref{eq2.23}) and (\ref{eq2.24}) and let $h:\M\to(0,h_0)$
be an exhaustion
function, satisfying the following conditions:
\smallskip

$a_1)$ there is $h'>0$ such that $h^{-1}([0,h'])$ is compact and $h$ is a 
solution of (\ref{eq2.19}) in the open set $h^{-1}((h',h_0));$
\smallskip

$a_2)$ for a.e. $t_1,\;t_2 \in (h',h_0)$, $t_1<t_2$,
$$
\int\limits_{\Sigma_h(t_2)}\langle{{\nabla h}\over {|\nabla h|}},
A(x,\nabla h)\rangle\, d{\cal H}^{n-1}
=\int\limits_{\Sigma_h(t_1)}\langle{{\nabla h}\over {|\nabla h|}}, 
A(x,\nabla h) \rangle\, d{\cal H}^{n-1}.
$$
\smallskip

\noindent
Here $d{\cal H}^{n-1}$ is the element of the $(n-1)-$dimensional 
Hausdorff measure on $\Sigma_h .$
Exhaustion functions with these properties will be called {\it the special
exhaustion functions of $\M$ with respect to $A$}. In most cases the
mapping $A$ will be the $p-$Laplace operator (\ref{ieq3.21}).

Since the unit vector $\nu={{\nabla h}/{|\nabla h|}}$ is
orthogonal to the
$h$--sphere $\Sigma_h$, the condition $a_2)$ means that
the flux of the vector field
 $A(m,\nabla h)$ through $h$--spheres $\Sigma_h(t)$ is constant.

{\it Suppose that the function $A(m,\xi)$ is continuously differentiable. If
\smallskip

$b_1)$ $h\in C^{2}({\cal M}\setminus K)$ and satisfies equation }
(\ref{eq2.19}), {\it and} 
\smallskip

$b_2)$ {\it at every point $m\in\M$ where $\partial\M$ has a tangent plane
$T_m(\partial\M)$ the condition
$$
\langle A(m,\nabla h(m)), \nu \rangle=0
$$
is satisfied where $\nu$ is a unit vector of the inner normal to the
boundary $\partial\M$,
then  $h$ is a special exhaustion function of the manifold }$\M.$
\smallskip

The proof of this statement is simple. Consider the open set 
$$
{\M}(t_1,t_2)=\{m\in {\M}:t_1<h(m)<t_2\},\quad 0< t_1< t_2 < \infty,
$$
with the boundary $\partial {\cal M} (t_1,t_2).$
%Using the Gauss formula, we have 
Using the Stokes formula, we have for noncritical values $t_1<t_2$ 
(for the definition of critical values of $C^k$-functions
see, for example, \cite[Part II, Chapter 2, \S 10]{DNF})
$$
\int\limits_{\Sigma_h(t_2)}\langle{{\nabla h}\over{|\nabla h|}},A(m,\nabla h)
\rangle d{\cal H}^{n-1} -
\int\limits_{\Sigma_h(t_1)}\langle {{\nabla h}\over{|\nabla h|}},A(m,\nabla h)
\rangle d{\cal H}^{n-1}
$$
$$
=\int\limits_{\partial{\cal M}(t_1,t_2)\cup\cup_{i=1,2}\Sigma_h (t_i)}
\langle \nu, A(m, \nabla h)\rangle d{\cal H}^{n-1} =
\int\limits_{\partial {\M}(t_1,t_2)}
\langle\nu, A(m,\nabla h)\rangle d{\cal H}^{n-1}
$$
$$= \int\limits_{{\M}(t_1,t_2)} {\rm div}\,A(m,\nabla h) *\1 = 0,
$$
and $a_2)$ follows.
 
\medskip

\begin{exmp}\label{3.26}
 We fix an integer $k$, $1\le k\le n,$ and set
$$
d_k(x)=\Bigl(\sum\limits_{i=1}^kx_i^2\Bigr)^{1/2}\,.
$$
It is easy to see that $|\nabla d_k(x)|=1$ everywhere in
${\bf R}^n\setminus\Sigma_0$ where $\Sigma_0=\{x\in {\bf R}^n:d_k(x)=0\}$. 
We shall call the set
$$B_k(t)=\{x\in {\bf R}^n:d_k(x)<t\}$$
a $k$-ball and the set
$$\Sigma_k(t)=\{x\in {\bf R}^n:d_k(x)=t\}$$
a $k$-sphere in ${\bf R}^n$.

We shall say that an unbounded domain $D\subset {\bf R}^n$ is $k$-admissible
if for each $t>\inf_{x\in D}d_k(x)$ the set $D\cap B_k(t)$ has compact closure.

It is clear that every unbounded domain $D\subset {\bf R}^n$ is
$n$-admissible. In the general case the domain $D$ is $k$-admissible
if and only if the function $d_k(x)$ is an exhaustion function of $D$.
It is not difficult to see that if a domain $D\subset {\bf R}^n$ is
$k$-admissible, then it is $l$-admissible for all $k<l<n$.

Fix $1\le k<n$. Let $\Delta$ be a bounded domain in the $(n-k)$-plane
$x_1=\ldots=x_k=0$ and let
$$D=\{x=(x_1,\ldots,x_k,x_{k+1},\ldots,x_n\in
{\bf R}^n:(x_{k+1},\ldots,x_n)\in\overline{\Delta}\}$$
be a closed domain in ${\bf R}^n$.

The domain $D$ is $k$-admissible. The  $k$-spheres $\Sigma_k(t)$ are
orthogonal to the boundary $\partial D$ and therefore
$\langle\nabla d_k,\nu\rangle=0$ everywhere on the boundary. The
function
$$
h(x)=\cases{\log d_k(x), &$p=k$,\cr
              d_k^{(p-n)/(p-1)}(x), &$p\neq k$,\cr}
$$
is a special exhaustion function of the domain $D$. Therefore for $p\ge
k$ the domain $D$ is of $p$-parabolic type and for $p<k$ $p$-hyperbolic
type.
\end{exmp}
\bigskip

\begin{exmp}
Fix $1\le k<n$. Let $\Delta$ be a bounded domain 
in the plane
$x_1=\ldots=x_k=0$ with a piecewise smooth
boundary and let
\eqb
\label{D}
D=\{x=(x_1,\ldots,x_n)\in
{\bf R}^n:(x_{k+1},\ldots,x_n)\in\Delta\} = {\bf R}^{n-k}\times\Delta
\eqe
be the cylinder domain with base $\Delta.$

The domain $D$ is $k$-admissible. The $k$-spheres $\Sigma_k(t)$ are
orthogonal to the boundary $\partial D$ and therefore
$\langle\nabla d_k,\nu\rangle=0$ everywhere on the boundary, where
$d_k$ is as in Example \ref{3.26}.

Let $h=\phi (d_k)$ where $\phi$ is a $C^2-$function. We have $\nabla
h=\phi^{\prime}\;\nabla d_k$ and
$$
\sum_{i=1}^n {{\partial}\over{\partial x_i}}\Bigl( |\nabla h|^{n-2}\;
{{\partial h}\over {\partial x_i}}\Bigr)=
\sum_{i=1}^k {{\partial}\over{\partial x_i}}\Bigl( (\phi^{\prime})^{n-1} \;
{{\partial d_k}\over{\partial x_i}}\Bigr)
$$
$$
=(n-1)\;(\phi^{\prime})^{n-2}\;\phi^{\prime\prime} + {{k-1}\over {d_k}}\;
(\phi^{\prime})^{n-1}.
$$
 From the equation
$$
(n-1)\;\phi^{\prime\prime} + {{k-1}\over {d_k}}\;\phi^{\prime}=0
$$
we conclude that  {\it the function}
\begin{equation}
\label{h}
h(x)= \bigl(d_k(x)\bigr)^{{n-k}\over{n-1}}
\end{equation}
{\it satisfies the equation} (\ref{ieq3.21}) {\it in $D \setminus K$
and thus it  is a special
exhaustion function of the domain} $D.$
%for the $n-$harmonic equation}
%$$A(x,h) = |h|^{n-2} h.$$
\end{exmp}

\medskip

\begin{exmp}
Let $(r,\theta)$, where $r\ge 0$, $\theta\in S^{n-1}(1)$, be the
spherical coordinates in $R^n$. Let $U\subset S^{n-1}(1)$, $\partial U\ne
\emptyset,$  be an arbitrary domain on the unit sphere $S^{n-1}(1)$.
We fix $0\le r_1<\infty$ and
consider the domain
\eqb
\label{-3}
D=\{(r,\theta)\in R^n: r_1<r<\infty,\;\theta\in U\}.
\eqe
As above it is easy to verify that the given domain is
$n$--admissible and {\it the function}
\eqb
\label{log}
h(|x|)=\log {{|x|}\over {r_1}}
\eqe
{\it is a special exhaustion function of the domain} $D$ for $p=n$.  
\end{exmp}
\medskip

\medskip

\begin{exmp}
Fix $1\le n\le p$. Let $x_1,x_2,\ldots,x_n$ be an orthonormal system of
coordinates in ${\bf R}^n,$ $1\le n< p$. Let $D\subset {\bf R}^n$ 
be an unbounded domain with piecewise smooth boundary and let
$\B$ be an $(p-n)$-dimensional compact Riemannian manifold with or without
boundary. We consider the manifold $\M= D\times\B$.

We denote by $x\in D$, $b\in\B$, and $(x,b)\in \M$ the points of the
corresponding manifolds. Let $\pi:D\times\B\to D$ and $\eta:D\times\B\to
\B$ be the natural projections of the manifold $\M$.

Assume now that the function $h$ is a function on the domain $D$ 
satisfying the conditions $b_1)$, $b_2)$ and the equation (\ref{ieq3.21}). We
consider the function $h^*=h\circ\pi:\M\to(0,\infty)$.

We have
$$
\nabla h^*=\nabla(h\circ\pi)=(\nabla_x h)\circ\pi
$$
and
$$
\diver(|\nabla h^*|^{p-2}\nabla h^*)=\diver\bigl(|\nabla(h\circ\pi)|^{p-2}
\nabla(h\circ\pi)\bigr)
$$
$$
=\diver\bigl(|\nabla_x h|^{p-2}\circ\pi(\nabla_x h)\circ\pi\bigr)
=\Bigl(\sum_{i=1}^n{\partial\over\partial x_i}\bigl(|\nabla_x h|^{p-2}
{\partial h\over\partial x_i}\bigr)\Bigr)\circ\pi.
$$
Because $h$ is a special exhaustion function of $D$ we have
$$
\diver(|\nabla h^*|^{p-2}\nabla h^*)=0.
$$

Let $(x,b)\in\partial\M$ be an arbitrary
point where the boundary $\partial\M$ has a tangent hyperplane and let $\nu$
be a unit normal vector to $\partial\M$.

If $x\in\partial D$, then $\nu=\nu_1+\nu_2$ where the vector 
$\nu_1\in {\bf R}^k$ is orthogonal to $\partial D$ and $\nu_2$ is 
a vector from $T_b(\B)$. Thus
$$
\langle\nabla h^*,\nu\rangle=\langle(\nabla_x h)\circ \pi,\nu_1\rangle=0,
$$
because $h$ is a special exhaustion function on $D$ and satisfies
the property $b_2)$ on $\partial D$.
If $b\in\partial\B$, then the vector $\nu$ is orthogonal to $\partial\B
\times {\bf R}^n$ and
$$
\langle\nabla h^*,\nu\rangle=\langle(\nabla_x h)\circ\pi,\nu\rangle=0,
$$
because the vector $(\nabla_x h)\circ\pi$ is parallel to ${\bf R}^n$.

The other requirements for a special exhaustion
function for the manifold $\M$ are easy to verify.

Therefore, {\it the function }
\eqb
\label{h^*}
h^*=h^*(x,b)=h\circ\pi:\M\to (0,\infty)
\eqe
{\it is a special exhaustion function on the manifold} $\M=D\times\B$.
\end{exmp}
\medskip

\begin{exmp}\label{3.25}
 Let $\A$ be a compact Riemannian manifold, $\dim
\A=k,$ with piecewise smooth boundary or without boundary. We consider the
Cartesian product $\M=\A\times {\bf R}^n$, $n\ge 1$. We denote by $a\in\A$,
$x\in {\bf R}^n$ and $(a,x)\in\M$ the points of the corresponding spaces. It
is easy to see that the function
$$
h(a,x)=\cases{\log |x|, &$p=n$,\cr
                  |x|^{p-n\over p-1}, &$p\neq n$,\cr}
$$
is a special exhaustion function for the manifold $\M$. Therefore, for
$p\ge n$ the given manifold is of $p$-parabolic type and for $p<n$
$p$-hyperbolic type.
\end{exmp}
\medskip

\begin{exmp}\label{3.26a}
Let $(r,\theta)$, where $r\ge 0$, $\theta\in S^{n-1}(1)$, be the
spherical coordinates in ${\bf R}^n$. Let $U\subset S^{n-1}(1)$ be an arbitrary
domain on the unit sphere $S^{n-1}(1)$. We fix $0\le r_1<r_2<\infty$ and
consider the domain
$$
D=\{(r,\theta)\in {\bf R}^n: r_1<r<r_2,\;\theta\in U\}
$$
with the metric
\eqb
\label{eq3.27a}
ds^2_\M=\alpha^2(r)dr^2+\beta^2(r) dl_{\theta}^2,
\eqe
where $\alpha (r),\,\beta (r)>0$ are $C^0$-functions on $[r_1,r_2)$ and
$dl_{\theta}$ is an element of length on $S^{n-1}(1)$.

The manifold $\M=(D,ds^2_\M)$ is a warped  Riemannian product. In the case 
$\alpha(r)\equiv 1$, $\beta(r)=1$, and $U=S^{n-1}$ the manifold
$\M$ is isometric to a cylinder in ${\bf R}^{n+1}$. In the case 
$\alpha (r)\equiv 1$, $\beta(r)=r$, $U=S^{n-1}$ the manifold 
$\M$ is a spherical annulus in ${\bf R}^{n}$.

The volume element in the metric (\ref{eq3.27a}) is given by the
expression
$$
d\sigma_\M=\alpha(r)\,\beta^{n-1}(r)\,dr\,dS^{n-1}(1).
$$
If $\phi(r,\theta)\in C^1(D)$, then the length of the gradient $\nabla\phi$
in $\M$ takes the form
$$
|\nabla\phi|^2={1\over\alpha^2}(\phi^{\prime}_r)^2+{1\over \beta^2}
|\nabla_{\theta}\phi|^2,
$$
where $\nabla_{\theta}\phi$ is the gradient in the metric of the unit
sphere $S^{n-1}(1)$.

For the special exhaustion function $h(r,\theta)\equiv h(r)$ the
equation (\ref{ieq3.21}) reduces to the following form
$$
{d\over dr}\left(\Bigl({1\over\alpha(r)}\Bigr)^{p-1}
\bigl(h^{\prime}_r(r)\bigr)^{p-1}\beta^{n-1}(r)\right)=0.
$$
Solutions of this equation are the functions
$$
h(r)=C_1\int\limits_{r_1}^r {\alpha(t)\over \beta^{n-1\over p-1}(t)}\,dt
+C_2
$$
where $C_1$ and $C_2$ are constants.

Because the function $h$ satisfies obviously the boundary condition
$a)_2$ as well as the other conditions of (\ref{3.20}), we see
that under the assumption
\eqb
\label{eq3.28a}
\int\limits^{r_2}{\alpha(t)\over \beta^{n-1\over p-1}(t)}\,dt=\infty
\eqe
the function
\eqb
\label{eq3.29a}
h(r)=\int\limits_{r_1}^{r}{\alpha(t)\over \beta^{n-1\over p-1}(t)}\,dt
\eqe
is a special exhaustion function on the manifold $\M$.
\end{exmp}
\medskip

\begin{thm}{}\label{3.23}
 Let $h:\M\to(0,h_0)$ be a special exhaustion function of
a boundary set $\xi$ of the manifold $\M$. Then

(i) if $h_0=\infty$, the set $\xi$ is of $p$-parabolic type,

(ii) if $h_0<\infty$, the set $\xi$ is of $p$-hyperbolic type.
\end{thm}

{\bf Proof.} Choose $0<t_1<t_2< h_0$ such that $K\subset B_h(t_1)$. 
We need to estimate the $p$-capacity of the condenser
$(B_h(t_1),\M\setminus B_h(t_2);\M)$. We have
\eqb
\hbox{\rm cap}_p(\overline{B}_h(t_1),\M\setminus B_h(t_2);\M)=
{J\over (t_2-t_1)^{p-1}}\label{eq3.24}
\eqe
where
$$J=\int\limits_{\Sigma_h(t)}|\nabla h|^{p-1}d{\cal H}_{\M}^{n-1}$$
is a quantity independent of $t>h(K)=\sup\{h(m):m\in K\}$.
Indeed, for the variational problem (\ref{eq3.6}) we choose the function
$\varphi_0$, $\varphi_0(m)=0$ for $m\in B_h(t_1)$,
$$\varphi_0(m)={h(m)-t_1\over t_2-t_1},\ m\in B_h(t_2)\setminus B_h(t_1)$$
and $\varphi_0(m)=1$ for $m\in\M\setminus B_h(t_2)$. Using the
Kronrod--Federer formula \cite[Theorem 3.2.22]{Fe}, we get
$$
\begin{array}{ll}
\hbox{\rm cap}_p(B_h(t_1),\M\setminus
B_h(t_2);\M)&\le\displaystyle\int\limits_{\M}|\nabla\varphi_0|^p*\1_{\M} \\ \\
\quad&
\le{1\over (t_2-t_1)^p}\displaystyle\int\limits_{t_1<h(m)<t_2}|\nabla
h(m)|^p*\1_{\M} \\ \\
\quad&=\displaystyle\int\limits_{t_1}^{t_2}dt\displaystyle
\int\limits_{\Sigma_h(t)}
|\nabla h(m)|^{p-1}d{\cal H}_{\M}^{n-1}.\\ \\
\end{array}
$$

Because the special exhaustion function satisfies the equation
(\ref{ieq3.21}) and the boundary condition $a)_2$, one obtains for arbitrary
$\tau_1,\tau_2$, $h(K)<\tau_1<\tau_2<h_0$
$$
\int\limits_{\Sigma_h(t_2)}|\nabla
h|^{p-1}d{\cal H}_{\M}^{n-1}-\int\limits_{\Sigma_h(t_1)}|\nabla h|^{p-1}
d{\cal H}_{\M}^{n-1}=
$$
$$
=\int\limits_{\Sigma_h(t_2)}|\nabla h|^{p-2}
\langle \nabla h, \nu \rangle d{\cal H}_{\M}^{n-1}
-\int\limits_{\Sigma_h(t_1)}|\nabla h|^{p-2}\langle \nabla h, 
\nu \rangle d{\cal H}_{\M}^{n-1}=
$$
$$
=\int\limits_{t_1<h(m)<t_2}\diver_{\M} (|\nabla h|^{p-2}\nabla h)*\1_{\M}=0.
$$
Thus we have established the inequality
$$\hbox{\rm cap}_p(B_h(t_1),\M\setminus B_h(t_2);\M)\le {J\over
(t_2-t_1)^{p-1}}\,.$$

By the conditions, imposed on the special exhaustion function, the
function $\varphi_0$ is an extremal in the variational problem (\ref{eq3.6}).
Such an extremal is unique and therefore the preceding inequality holds in
fact with equality. This conclusion proves the equation (\ref{eq3.24}).

If $h_0=\infty$, then letting $t_2\to\infty$ in (\ref{eq3.24}) we conclude the
parabolicity of the type of $\xi$. Let $h_0<\infty$. Consider an
exhaustion $\{\U_k\}$ and choose $t_0>0$ such that the $h$-ball $B_h(t_0)$ 
contains the compact set $K$.

Set $t_k=\sup\nolimits_{m\in\partial\U_k}h(m)$.
Then for $t_k>t_0$ we have
$$
\hbox{\rm cap}_p(\overline{U}_{k_0},\U_k;\M)\ge 
\hbox{\rm cap}_p(B_h(t_0),B_h(t_k);\M)= J/(t_k-t_0)^{p-1}\,,
$$
and hence
$$\liminf_{k\to\infty} \hbox{\rm cap}_p(\overline{U}_{k_0},\U_k;\M)\ge
J/(h_0-t_0)^{p-1}>0,$$
and the boundary set $\xi$ is of $p$-hyperbolic type.
$\Box$

\bigskip

%%%SECTION
%%%SECTION
%%%SECTION

\cc
\section{Energy integral}{}\label{sec4}
The fundamental result of this section is an estimate for the rate of
growth of the energy integral of forms of the class $\WT_2$ on noncompact
manifolds under various boundary conditions for the forms. As an
application we get Phragm\'en--Lindel\"of type theorems for the forms of
this class and we prove some generalizations of the classical theorem of
Ahlfors concerning the number of distinct asymptotic tracts of an entire
function of finite order.
\bigskip

\begin{nonsec}{Boundary conditions }\label{4.1}
\end{nonsec}
 Let $\M$ be an $n$-dimensional
Riemannian manifold with nonempty boundary $\partial\M$. We will fix a
closed differential form $w$, $\deg w=k$, $1\le k\le n$, $w\in 
L_{\loc}^p(\M)$ of class $\WT_1$ and the complementary closed form
$\theta$, $\deg\theta=n-k$, $\theta\in L_\loc^q(\M)$, satisfying
the condition (\ref{eq2.4}). We assume that there exists a 
differential form $Z\in W_{p,\loc}^1$ with continuous coefficients
for which $dZ=w$.
\medskip

Let $h:\M\to(0,h_0)$ be an exhaustion function of $\M$. As before we let
$B_h(\tau)$ be an $h$-ball and $\Sigma_h(\tau)$ an $h$-sphere.
\bigskip

\begin{nonsec}{Dirichlet condition with zero boundary values }\label{4.2}
\end{nonsec}
We shall say that the form $Z\in W_{p,\loc}^1$
\footnote {with continuous coefficients and such that $dZ=w$} 
satisfies Dirichlet's condition with zero boundary values on 
$\partial \M$ if for every differential form $v\in L_{\loc}^q (\M)$, 
$\deg v=n-k$, and for almost every $\tau\in(0,h_0)$
\eqb
\label{eq4.3}
\int\limits_{B_h(\tau)} w\wedge v+(-1)^{k-1}\int\limits_{B_h(\tau)} Z\wedge
dv=\int\limits_{\Sigma_h(\tau)}Z\wedge v.
\eqe

In particular, the form $Z$ satisfies the boundary condition (\ref{eq4.3}),
if its coefficients are continuous and if its support does not intersect
with $\partial\M$, that is
\eqb
\label{eq4.4}
\supp Z\cap\partial\M=\emptyset,\ \hbox{\rm where}\
\supp Z=\overline{\{m\in\M:Z(m)\neq 0\}}.
\eqe

If $\M$ is compact then (\ref{eq4.3}) takes the form
\eqb
\label{eq4.3cm}
\int\limits_{\M} w\wedge v+(-1)^{k-1}\int\limits_{\M} Z\wedge
dv=0.
\eqe
\bigskip

\begin{nonsec}{Neumann condition with zero boundary values }\label{4.5}
\end{nonsec}
We shall say that a form $Z$ satisfies Neumann's condition with zero boundary
values, if for every differential form $v\in W_{p,\loc}^1(\M)$, $\deg
v=k-1$, and for almost every $\tau\in(0,h_0)$
\eqb
\label{eq4.6}
\int\limits_{B_h(\tau)} dv\wedge\theta=\int\limits_{\Sigma_h(\tau)}
v\wedge\theta.
\eqe

If $\M$ is compact then (\ref{eq4.6}) takes the form
\eqb
\label{eq4.6cm}
\int\limits_{\M} dv\wedge\theta=0.
\eqe
\bigskip

\begin{nonsec}{Mixed zero boundary condition }\label{4.7a}
\end{nonsec}
We shall say that a form $Z$ satisfies mixed zero boundary condition if
for an arbitrary function $\phi\in C^1(\M)$ and for almost every
$\tau\in (0,h_0)$ we have
\eqb
\label{eq4.8a}
\int\limits_{B_h(\tau)}\phi w\wedge\theta+(-1)^{n-1}\int\limits_{B_h(\tau)}
Z\wedge\theta\wedge d\phi=\int\limits_{\Sigma_h(\tau)}\phi Z\wedge\theta.
\eqe

If $\M$ is compact then (\ref{eq4.8a}) takes the form
\eqb
\label{eq4.8acm}
\int\limits_{\M}\phi w\wedge\theta+(-1)^{n-1}\int\limits_{\M}
Z\wedge\theta\wedge d\phi=0.
\eqe
 \smallskip

%\begin{subsec}{}\label{4.7}
%\end{subsec}
 We assume that the form
\eqb
\label{eq4.8}
Z\in C^2(\hbox{\rm int}\M)\cap C^1(\partial\M)
\eqe
has the property (\ref{eq4.3}). On the basis of Stokes' formula 
(the standard Stokes formula with generalized derivatives) we conclude
that for almost every $\tau\in(0,h_0)$
$$
\int\limits_{B_h(\tau)} dZ\wedge v+(-1)^{k-1}\int\limits_{B_h(\tau)} Z\wedge
dv=\int\limits_{\partial B_h(\tau)} Z\wedge v
$$
holds. Therefore we get
$$
\int\limits_{\partial B_h(\tau)\setminus\Sigma_h(\tau)} Z\wedge v=0\quad
\hbox{for all}\;\; v\in W_{q,\loc}^1(\M).
$$
This implies that the restriction of $Z$ onto the boundary $\partial\M$ is the
zero form, i.e.
\eqb
\label{eq4.9}
Z\,|_{\partial\M}(m)=0\quad\hbox{at every point}\;\;m\in\partial\M.
\eqe

\smallskip

%\begin{subsec}{}\label{4.10}
%\end{subsec}
We next clarify the geometric meaning of the condition (\ref{eq4.9}). We
assume that $m\in\partial\M$ is a point where the boundary $\partial\M$
has a tangent plane $T_m(\partial\M)$ and that the form $Z$ satisfies
the regularity condition (\ref{eq4.8}) in some neighborhood of the point $m$.
\smallskip

\begin{prop}{}\label{4.11} If a form $Z$ is simple at a point $m\in\M$,
then the condition (\ref{eq4.9}) is fulfilled if and only if the form $Z$
is of the form
\eqb
Z=\omega\wedge dx_n,\label{eq4.12}
\eqe
where $\omega$ is a form, ${\rm deg}\,\omega={\rm deg}\,Z -1$.
\end{prop}
\smallskip

{\bf Proof.} 
We give an orthonormal system of coordinates $x_1,\ldots,x_n$ at the
point $m$ such that the hyperplane $T_m(\partial\M)$ is given by the
equation $x_n=0$. Let ${\rm deg}\,Z=l.$
Because the form $Z$ is simple, we can represent it as follows
\eqb
\label{eq4.13}
\begin{array}{ll}
Z&=\Bigl(\sum_{i=1}^l a_{1,i}dx_i+a_{1,n}dx_n\Bigr)\wedge\ldots\\ \\
\quad &\ldots\wedge
\Bigl(\sum_{i=1}^l
a_{l,i}dx_i + a_{l,n}dx_n\Bigr),\\ \\
\end{array}
\eqe
where $a_{i,j}=a_{i,j}(m)$ are some constants. The condition (\ref{eq4.9}) can
now be rewritten as follows
$$
\Bigl(\sum_{i=1}^l a_{1,i}dx_i\Bigr)\wedge\ldots\wedge
\Bigl(\sum_{i=1}^l a_{l,i}dx_i\Bigr)=0,
$$
and we easily obtain (\ref{eq4.12}).

The proof of the converse implication is obvious. $\Box$
\bigskip

%\begin{subsec}{}\label{4.14}
%\end{subsec}
We next clarify the geometric meaning of the Neumann condition
(\ref{eq4.6}). We fix the forms
$$
Z,v\in C^2(\hbox{\rm int}\M)\cap C^1(\partial\M).
$$
By Stokes' formula we have for almost every $ \tau \in (0,h_0)$
$$\int\limits_{\partial B_h(\tau)}v\wedge\theta=
\int\limits_{B_h(\tau)} dv\wedge\theta
+(-1)^{k-1} \int\limits_{B_h(\tau)} v\wedge d\theta.$$
Because the form $\theta$ is closed, condition (\ref{eq4.6}) gives
$$\int\limits_{\Sigma_h(\tau)} v\wedge\theta=0\quad\mbox{for all}\quad v\in
W_{p,\loc}^1\,.$$
Therefore
\eqb
\theta\,|_{\partial\M}(m)=0\label{eq4.15}
\eqe
at every point $m\in\partial\M$. 
\smallskip

Exactly in the same way we verify that the mixed zero boundary
condition (\ref{eq4.8a}) is equivalent to the condition
\eqb
Z \wedge \theta\,|_{\partial\M}(m)=0 \label{eq4.15b}
\eqe
at every point  $m\in\partial\M$. 
\bigskip

%\begin{subsec}{}\label{4.16}
%\end{subsec}
Consider the case of quasilinear equations (\ref{eq2.19}).
Let $m\in\partial\M$ be a regular point and let $x_1,\ldots,x_n$ be local
coordinates in a neighborhood of this point.  We have
\eqb\label{4.16}
\begin{array}{ll}
\theta&=*\sum_{i=1}^nA_i(m,\nabla f(m))dx_i \\ \\
\quad&=\sum_{i=1}^n (-1)^{i-1} A_i(m,\nabla f(m))dx_1\wedge\ldots
\wedge\widehat{dx_i}\wedge\ldots\wedge dx_n.\\ \\ 
\end{array}
\eqe
We set $Z=f$. In the case (\ref{eq4.3}) we choose 
$v=\phi\theta$ where $\theta$ is an arbitrary locally Lipschitz 
function. We obtain
$$
\int\limits_{B_h(\tau)} \phi df\wedge\theta+ \int\limits_{B_h(\tau)}
fd(\phi\theta)=\int\limits_{\Sigma_h(\tau)} \phi f\theta
$$
and further
\eqb
\label{eq4.17}
\int\limits_{B_h(\tau)}\sum_{i=1}^n(\phi f)_{x_i}A_i(m,\nabla f)*\1
=\int\limits_{\Sigma_h(\tau)}\phi f\theta,\quad
\hbox{for all}\;\phi.
\eqe
This condition characterizes generalized solutions of the equation
(\ref{eq2.19}) with zero Dirichlet boundary condition on $\partial\M$.

On the other side, choosing in the case of the Neumann condition (\ref{eq4.6})
for $v$ an arbitrary locally Lipschitz function $\phi$ we get for almost
every $\tau\in(0,h)$
\eqb
\label{eq4.18}
\int\limits_{B_h(\tau)}\langle\nabla\phi, A(m,\nabla f)\rangle*\1=
\int\limits_{\Sigma_h(\tau)}\phi\langle A(m, \nabla f),
 \nu\rangle d{\cal H}_{\M}^{n-1}
\eqe
which characterizes generalized solutions of the equation (\ref{eq2.19}) with
zero Neumann boundary conditions on $\partial\M$.

It is easy to see that at every point of the boundary we have
$$
(-1)^{i-1}dx_1\wedge\ldots\wedge\widehat{dx_i}\wedge\ldots\wedge dx_n
\,|_{\partial\M}=\cos(\nu,x_i)\;d{\cal H}^{n-1}_{\M},
$$
where $(\nu,x_i)$ is the angle between the inner normal vector $\nu$ to
$\partial\M$ and the direction $0x_i$; $d{\cal H}_{\M}^{n-1}$ is the 
element of surface area on $\M$..

Thus, at a regular boundary point, the condition (\ref{eq4.15}) 
is equivalent to the requirement
$$
\langle A(m,\nabla f(m)),\nu\rangle=0.
$$

Using (\ref{eq4.9}) we see that the condition (\ref{eq4.8a}) is
equivalent to the traditional mixed boundary condition at regular
boundary points.
\bigskip

\begin{nonsec}{Maximum principle for $\WT$-forms }\label{4.19}
\end{nonsec}
 Let $\M$ be a
compact Riemannian manifold with nonempty boundary, $\dim\M=n$, and let
$v\in L_{\loc}^p$, $\deg w=k$, $1\le k\le n$, be a 
differential form of class $\WT_1$ on $\M$. Let $\theta$,
$\deg\theta=n-k$, be a form complementary to the form $w$.
\smallskip

\begin{thm}{}\label{4.20}
Suppose that there exists a differential form
$Z\in W_{p,\loc}^1(\M)$, $dZ=w$ on $\M$. If either (\ref{eq4.3cm})
or (\ref{eq4.6cm}) holds, then $\theta\equiv 0$ on $\M$.
\end{thm}
\smallskip

{\bf Proof.} We assume that (\ref{eq4.3cm}) holds and set $v=\theta$. Then
(\ref{eq4.3cm}) yields
$$
\int\limits_{\M}w\wedge\theta=0.
$$
Because
$$
\begin{array}{ll}
(-1)^{k(n-k)}(-1)^{k(n-k)}*(w\wedge\theta)&=(-1)^{k(n-k)}*^{-1}
(w\wedge\theta) \\ \\
\quad&=*^{-1}(w\wedge (-1)^{k(n-k)}\theta) \\ \\
\quad &=*^{-1}(w\wedge
*(*\theta))=\langle w ,*\theta\rangle,\\ \\
\end{array}
$$
we get
$$\int\limits_{\M}w\wedge\theta =\int\limits_{\M} *(w\wedge \theta)*\1
=\int\limits_{\M}\langle w,*\theta\rangle *\1.$$
Using (\ref{eq2.5}) we deduce
\eqb
0=\int\limits_{\M} w\wedge\theta\ge {\nu}_0\int\limits_{\M}
|\theta|^q*\1.\label{eq4.21}
\eqe

We assume that the boundary condition (\ref{eq4.6cm}) holds. Choose $v=Z$. Then
(\ref{eq4.6cm}) gives
$$\int\limits_{\M}w \wedge\theta=0.$$
As above, we arrive at the inequality (\ref{eq4.21}). This inequality implies
that $\theta\equiv 0$ on $\M$. $\Box$

\smallskip

 In order to illustrate Theorem \ref{4.20} we consider the example
of generalized solutions $f\in W_{p,\loc}^1(\M)$ of the equation
(\ref{eq2.19}) under the condition 
\eqb
\nu_0\,|A(m,\xi)|^p\le \langle\xi, A(m,\xi)\rangle \label{eq2.18}
\eqe
for all $\xi\in T_m(\M)$
with the constants $p>1$ and $\nu_0>0$.

Setting $Z=f$ we get
\smallskip

\begin{cor}\label{4.23}
 Suppose that the manifold $\M$ is compact and
the boundary $\partial\M$ is not empty. If the function $f$ satisfies
the condition (\ref{eq4.17}) or (\ref{eq4.18}), then 
$f\equiv {\rm const}$ on $\M$.
\end{cor}

\bigskip

\bigskip

%%%SECTION
%%%SECTION
%%%SECTION

\cc
\section{Estimates for energy integral. Applications}{}\label{4.31}
This Chapter is devoted to Phragm\'en-Lindel\"of and Ahlfors theorems 
for differential forms.
\smallskip

\begin{nonsec}{Basic theorem }{}\label{subsbas}
\end{nonsec}
Let $\M$ be a noncompact Riemannian manifold, $\dim\M=n$. We consider a
class $\F$ of differential forms $Z\in W_{p,\loc}^1(\M)$, $\deg
Z=k-1$, such that the form $dZ=w$ satisfies the conditions (\ref{eq2.2})
and is in the class $\WT_2$. Let $\theta\in L_{\loc}^q$ be a form
satisfying the condition (\ref{eq2.4}), complementary to $w$.

If the boundary $\partial\M$ is nonempty then we shall assume that the
form $Z$ satisfies on $\partial\M$ some boundary condition $B$. In the
case considered below such a boundary condition can be any of the
conditions (\ref{eq4.3}), (\ref{eq4.4}), (\ref{eq4.6}), (\ref{eq4.8a}).
We shall denote by $\F_B(\M)$ the set of forms $Z$, $dZ\in\WT_2$,
satisfying the boundary condition $B$ on $\M$. In particular, below we
shall operate with the classes $\F_D$, $\F_0$, $\F_{N}$, and $\F_{DN}$
forms corresponding to the boundary conditions (\ref{eq4.3}), (\ref{eq4.4}),
(\ref{eq4.6}), (\ref{eq4.8a}), respectively.

We fix a locally Lipschitz exhaustion function $h:\M\to(0,h_0)$,
$0<h_0\le\infty$. Let
$\tau\in (0,h_0)$ and let $B_h(\tau)$ be an $h$-ball, and
$\Sigma_h(\tau)$ its boundary sphere as before.

We introduce a characteristic $\epsilon(\tau)$ setting
\eqb
\label{eq4.34}
\epsilon(\tau;\F_B)=
\inf{\displaystyle\int\limits_{\Sigma_h(\tau)}|w|^p
|\nabla h|^{-1}d{\cal H}_{\M}^{n-1}
\over \Bigl|\displaystyle\int\limits_{\Sigma_h(\tau)}\langle
Z,*\theta\rangle\, d{\cal H}_{\M}^{n-1}\Bigr|}
\eqe
where the infimum is taken over all $Z\in\F_B(\M)$, $Z\neq 0$.

Some estimates of (\ref{eq4.34}) are given in \cite{klv} and \cite{mmv7}.

Under these circumstances we have
\medskip

\begin{thm}{}\label{4.35}
Suppose that the form $Z \in \F_B(\M)$ satisfies one of the boundary
condition (\ref{eq4.3}), (\ref{eq4.6}), or (\ref{eq4.8a}).
Then for almost all $\tau\in(0,h_0)$ and for an arbitrary $\tau_0$ the
following relation holds
\eqb
\label{eq4.36}
{d\over d \tau}\Bigl(I(\tau)\exp\Bigl\{-\nu_1\int
\limits_{\tau_0}^{\tau}\epsilon(t;\F_B)\,dt\Bigr\}\Bigr)\ge 0,
\eqe
where
$$I(\tau)=\int\limits_{B_h(\tau)}|w|^p*\1.$$

In particular, for all $\tau_1<\tau_2$ we have
\eqb
I(\tau_1)\le I(\tau_2)\exp\{-{\nu_1}\int\limits_{\tau_1}^{\tau_2}\epsilon(t)\, dt\}.\label{eq4.37}
\eqe
\end{thm}
\smallskip

{\bf Proof.} The Kronrod--Federer formula yields
\eqb
I(\tau)=\int\limits_0^{\tau}dt\int\limits_{\Sigma_h(\tau)}|w|^p{d{\cal H}_{\M}^{n-1}\over|
\nabla h|}\label{eqKFf}
\eqe
and, in particular, the function $I(\tau)$ is absolutely continuous on closed
intervals of $(0,h_0)$. Now it is  enough to prove the inequality
\eqb
{d\over d\tau}I(\tau)\ge \nu_1\,I(\tau)\,\epsilon(\tau).\label{eq4.38}
\eqe
 From (\ref{eqKFf}) we have for almost every $\tau\in(0,h_0)$ 
\eqb
{d\over d\tau}I(\tau)=\int\limits_{\Sigma_h(\tau)}|w|^p{d{\cal H}_{\M}^{n-1}\over|\nabla
h|}.\label{eq4.39}
\eqe

By (\ref{eq2.7}) we obtain
$$
\begin{array}{ll}
I(\tau)&=\displaystyle\int\limits_{B_h(\tau)}|w|^p*\1\le
{\nu}_1^{-1}\displaystyle\int\limits_{B_h(\tau)}\langle w,*\theta\rangle*\1 \\  \\
\quad&={\nu}_1^{-1}\displaystyle\int\limits_{B_h(\tau)}w\wedge\theta=
{\nu}_1^{-1}\displaystyle\int\limits_{B_h(\tau)}
dZ\wedge\theta.
\end{array}
$$
However, the form $Z$ is weakly closed and satisfies one of the conditions (\ref{eq4.3}),
(\ref{eq4.4}), or (\ref{eq4.6}).
Therefore for a.e. $\tau\in (0,h_0)$,
$$
\int\limits_{B_h(\tau)}dZ\wedge\theta=\int\limits_{\Sigma_h(\tau)}
Z\wedge\theta.
$$
Thus we get
$$
I(\tau)\le {\nu}_1^{-1}\int\limits_{\Sigma_h(\tau)}\langle Z,*\theta\rangle\,d{\cal H}_{\M}^{n-1}.
$$

Further from (\ref{eq4.34}) it follows that
$$
\int\limits_{\Sigma_h(\tau)}|w|^p{d{\cal H}_{\M}^{n-1}\over|\nabla h|}\ge\epsilon
(\tau;\F_B)\Bigl|\int\limits_{\Sigma_h(\tau)}\langle
Z,*\theta\rangle\,d{\cal H}_{\M}^{n-1}\Bigr|.
$$

Combining the above inequalities we obtain
$$
I(\tau)\le {{\nu_1^{-1}}\over \epsilon(\tau; \F_B)}
\int\limits_{\Sigma_h(\tau)}|w|^p{d{\cal H}_{\M}^{n-1}\over|\nabla h|}.
$$
This inequality together with the equality (\ref{eq4.39}) yields
$$
I(\tau)\le{{\nu_1^{-1}} \over \epsilon(\tau; \F_B)}
{d\over d\tau}I(\tau).
$$
We thus obtain the desired conclusion (\ref{eq4.38}).
$\Box$

\smallskip

We shall need also some other estimates of the energy integral. We now
prove the first of these inequalities. Denote by $\F(B_h(\tau))$
the set of all differential forms
\eqb
\label{eq4.39a}
Z_0 \in C^1(B_h(\tau)),\quad\deg Z_0 = k-1,\quad dZ_0 =0,
\eqe
such that for almost every $\tau \in (0,h_0)$ and for an arbitrary
Lipschitz function $\phi$ the following formula holds
\eqb
\label{eq4.39b}
\int\limits_{\Sigma_h(\tau)}\phi\, Z_0\wedge\theta=\int_{B_h(\tau) } d\phi
\wedge Z_0 \wedge\theta.
\eqe
\smallskip

\begin{thm}{}\label{4.40}
If the differential form $Z\in\F_B(\M)$, $dZ\in\WT_2$,
satisfies the boundary condition (\ref{eq4.3}), (\ref{eq4.6}), or
(\ref{eq4.8a}), then for all $\tau_1<\tau_2<h_0$ and for an arbitrary
form $Z_0\in\F(B_h(\tau_2))$ the following relation holds
\eqb
\label{eq4.41}
\nu_1\int\limits_{B_h(\tau_1)}|dZ|^p*\1\le {{p}\over \tau_2-\tau_1}
\int\limits_{B_h(\tau_2)\setminus B_h(\tau_1)} |\nabla h||(Z-Z_0)
\wedge\theta|*\1.
\eqe
\end{thm}
\smallskip

{\bf Proof.} We consider the function
$$
\phi(m)=\cases{\qquad 1 &for $m\in B_h(\tau_1)$,\cr
              \db{\tau_2-h(m)\over \tau_2-\tau_1}\de &for $m\in
                 B_h(\tau_2)\setminus B_h(\tau_1)$,\cr
               \qquad 0 &for $m\in\M\setminus B_h(\tau_2)$.\cr}
$$
Suppose that the form $Z$ satisfies the condition (\ref{eq4.3}). Setting
in (\ref{eq4.3}) $v=(\phi)^p Z\wedge\theta$ we get
$$
\int\limits_{B_h(\tau_2)}(\phi)^p w\wedge\theta+(-1)^{k-1} 
\int\limits_{B_h(\tau_2)} Z\wedge d(\phi)^p\wedge\theta=0,
$$
or
$$
\int\limits_{B_h (\tau_2)}(\phi)^p\langle w,*\theta\rangle\,*\1=(-1)^k \, p \,
\int\limits_{B_h(\tau_2)}(\phi)^{p-1}\,Z\wedge d\phi\wedge\theta.
$$

The function $(\phi)^p$ is locally Lipschitz on 
$\overline B_h(\tau_2)$ and $\phi|_{\Sigma_h(\tau_2) }=0$.
%Therefore on the basis of (\ref{eq4.39b}) we have
%$$
%\int\limits_{\partial\M\cap B_h(\tau_2)}(\phi)^p\,Z_0\wedge\theta=0.
%$$
Thus by (\ref{eq4.39b}) we get
$$
\int\limits_{B_h(\tau_2)}d\phi^p\wedge Z_0\wedge\theta=
\int\limits_{B_h(\tau_2)} d(\phi)^p\wedge Z_0 \wedge \theta+
$$
$$
+\int\limits_{B_h(\tau_2)}(\phi)^p\,dZ_0\wedge\theta+
\int\limits_{B_h(\tau_2)}(\phi)^p\,Z_0\wedge d\theta=
$$
$$
=\int\limits_{\Sigma_h(\tau_2)}(\phi)^p\,Z_0\wedge\theta=0.
$$

Hence we arrive at the relation
$$
\int\limits_{B_h(\tau_2)}(\phi)^p\langle w,*\theta\rangle\,*\1=
(-1)^k \, p \, \int\limits_{B_h(\tau_2)}(\phi)^{p-1}(Z-Z_0)\wedge
d\phi\wedge\theta,
$$
which by (\ref{eq2.7}) yields
\eqb
\label{eq4.42}
\nu_1\int\limits_{B_h(\tau_2)}(\phi)^p |w|^p *\1\le{{p}\over
\tau_2-\tau_1}\int\limits_{B_h(\tau_2)}(\phi)^{p-1}\,|(Z - Z_0)\wedge
\theta|\,|\nabla h|*\1.
\eqe

Observing that $\phi(m)=1$ for $m\in B_h(\tau_1)$ and $\phi(m)=0$ for
$m\in\M\setminus B_h(\tau_2)$ we obtain
$$
\nu_1\int\limits_{B_h(\tau_1)}|w|^p*\1\le{{p}\over \tau_2-\tau_1}\int
\limits_{B_h(\tau_2)\setminus B_h(\tau_1)}(\phi)^{p-1}|\nabla h|\,
|(Z-Z_0)\wedge\theta|*\1.
$$
Because $|\phi|\le 1$, the inequality (\ref{eq4.41}) follows.

Let the form $Z$ satisfy the condition (\ref{eq4.6}). We choose $v=(\phi)^p Z$
and observe that
$$
v|_{\Sigma_h(\tau_2)}=0.
$$
Then we get
$$
\int\limits_{B_h(\tau_2)}(\phi)^p\,w\wedge\theta=-\int\limits_{B_h(\tau_2)}
d(\phi)^p\wedge Z\wedge\theta=(-1)^k \, p \, \int\limits_{B_h(\tau_2)}
(\phi)^{p-1} Z\wedge d\phi\wedge\theta.
$$
Further details of the proof in this case are similar 
to those carried out above.

We assume that the form $Z$ satisfies the mixed boundary condition
(\ref{eq4.8a}). Observing that
$$
(\phi)^p|_{\Sigma_h(\tau_2)}=0,
$$
we get
$$
\int\limits_{B_h(\tau_2)}(\phi)^p\,w\wedge\theta=(-1)^n\int
\limits_{B_h(\tau_2)} Z\wedge\theta\wedge d(\phi)^p=(-1)^{n-k} \, p\, 
\int \limits_{B_h(\tau_2)}(\phi)^{p-1}\,Z\wedge d\phi\wedge\theta.
$$
Arguing as above we complete the proof of the theorem.
$\Box$

\smallskip

There is also an estimate for the energy integral which does not use
the complementary form $\theta$ of $dZ=w$. Such
an estimate is given in the next theorem.
\smallskip

\begin{thm}{}\label{4.43}
If the form $Z$, $dZ\in\WT_2(\M)$, satisfies on $\partial\M$ one of the
boundary conditions (\ref{eq4.3}), (\ref{eq4.6}), or (\ref{eq4.8a}), then
for all $0<\tau_1<\tau_2 < h_0$ and for an arbitrary form 
$Z_0\in\F(B_h(\tau_2))$ we have
\eqb
\label{eq4.43b}
\int\limits_{B_h(\tau_1)} |dZ|^p *\1\le
\Bigl({{p\nu_2}\over{(\tau_2-\tau_1)\nu_1}}\Bigr)^p \int\limits_{B_h(\tau_2)
\setminus B_h(\tau_1)}|\nabla h|^p|Z-Z_0|^p*\1.
\eqe
\end{thm}
\smallskip

{\bf Proof.} We use the earlier established relation (\ref{eq4.42}). 
We estimate the integral on the right hand side of (\ref{eq4.42}). By
(\ref{eq2.8}) we get
$$
\int\limits_{\M}(\phi)^{p-1}|\nabla h||(Z-Z_0)\wedge\theta*\1 \le
\int\limits_{\M}(\phi)^{p-1}|\nabla h||Z-Z_0||\theta|*\1
$$
$$
\le {\nu}_2\int\limits_{\M}(\phi)^{p-1}|\nabla h||Z-Z_0|
|w|^{p-1}*\1
$$
$$
\le {\nu}_2\Bigl(\int\limits_{\M} |\nabla h|^p|Z-Z_0|^p*\1
\Bigr)^{1/p} \Bigl(\int\limits_{\M}\phi^p|w|^p*\1\Bigr)^{(p-1)/p}.
$$
 From (\ref{eq4.42}) we get
$$
\Bigl({\nu_1\over\nu_2}\Bigr)^p \int\limits_{\M}|w|^p*\1 \le \Bigl(
{p\over \tau_2-\tau_1}\Bigr)^p \int\limits_{\M}|\nabla h|^p|Z-Z_0|^p*\1.
$$
Using the facts that $\phi=1$ on $B_h(\tau_1)$ and $\phi=0$ on
$\M\setminus B_h(\tau_2)$ we easily verify come to (\ref{eq4.43b}).
$\Box$

\bigskip

\begin{nonsec}{Phragm\'en--Lindel\"of\ \ theorem }\label{4.45}
\end{nonsec}
Let $\M$ be an $n$-dimensional noncompact Riemannian manifold with or
without boundary and let $w\in\WT_2$ be a differential form
as in (\ref{eq2.2}), $\deg w=k$, and $\theta$ its complementary form
as in (\ref{eq2.4}).

We assume that there exists a differential form $Z\in
W_{p,\loc}^1$ with $dZ=w$. If the boundary $\partial\M$ is
nonempty, then we shall assume that $Z$ satisfies the boundary condition
of Dirichlet (\ref{eq4.3}), Neumann's condition (\ref{eq4.6}), or the condition (\ref{eq4.8a}).

We fix a locally Lipschitz exhaustion function $h:\M\to(0,h_0)$,
$0<h_0\le\infty$. Let,
as above
$$
I(\tau;Z)=\int\limits_{B_h(\tau)}|d Z|^p\,*\1
$$
and let
$$
\mu(\tau;Z)=\inf \int\limits_{\tau<h(m)<\tau+1}|\nabla h|\,|(Z-Z_0)\wedge
\theta| *\1,
$$
$$
m(\tau;Z)=\inf\int\limits_{\tau<h(m)<\tau+1}|\nabla h|^p\,|Z-Z_0|^p *\1,
$$
where the infimum is taken over all closed forms $Z_0$, satisfying
conditions (\ref{eq4.39a}), (\ref{eq4.39b}) on $B_h(\tau)$.

The following theorem exhibits a generalization of the classical
Phragm\'en--Lindel\"of principle for holomorphic functions.
\smallskip

\begin{thm}{}\label{4.46}
Suppose that the form $Z$, $dZ\in\WT_2(\M)$, satisfies one of the
boundary conditions (\ref{eq4.3}), (\ref{eq4.6}) or
(\ref{eq4.8a}). The following alternatives hold: either the form $dZ=0$
a.e. on the manifold $\M$, or for all $\tau_0\in(0,h_0)$ we have
\eqb
\label{eq4.47a}
\liminf_{\tau\to h_0} I(\tau;Z)\,\exp\Bigl\{-\nu_1\int
\limits_{\tau_0}^{\tau}\epsilon(t;\F_B)\,dt\Bigr\}>0;
\eqe

\eqb
\label{eq4.47b}
\liminf_{\tau\to h_0} \mu(\tau;Z)\,\exp\Bigl\{-\nu_1\int\limits_{\tau_0}^{\tau}\epsilon (t;\F_B)\,dt\Bigr\}>0,
\eqe

\eqb
\label{eq4.47c}
\liminf_{\tau\to h_0} m(\tau;Z)\,\exp\Bigl\{-\nu_1\int\limits_{\tau_0}^{\tau}\epsilon (t;\F_B)\,dt\Bigr\}>0.
\eqe
\end{thm}
\smallskip

{\bf Proof.}
The property (\ref{eq4.47a}) follows readily from (\ref{eq4.37}) and is
presented here only for the sake of completeness.

By (\ref{eq4.37}) and (\ref{eq4.41}), for a.e. $\tau\in (\tau_0,h_0)$ we have
$$
I(\tau_0)\le I(\tau)\,\exp\Bigl\{-\nu_1\int\limits_{\tau_0}^{\tau}
\epsilon(t)\,dt\Bigr\}\le
$$
$$
\le p\,\nu_1^{-1}\int\limits_{B_h(\tau+1)\setminus B_h(\tau)}|\nabla h|\,
|(Z-Z_0)\wedge\theta|*\1\;\;\;\exp\Bigl\{-\nu_1
\int\limits_{\tau_0}^{\tau}\epsilon(t)\,dt\Bigr\}.
$$
Therefore we get
$$
I(\tau_0)\le p\,\nu_1^{-1}\,\mu(\tau;Z)\,\exp\Bigl\{-\nu_1\int\limits_{\tau_0}^{\tau} \epsilon(t)\,dt\Bigr\}.
$$

Analogously, using (\ref{eq4.43b}) we get
$$
I(\tau_0)\le\Bigl({\nu_2\over\nu_1}\Bigr)^p\,p^p\,m(\tau;Z)\,\exp\Bigl\{
-\nu_1\int\limits_{\tau_0}^{\tau}\epsilon(t)\,dt\Bigr\}.
$$

If we now assume that the form $w\not\equiv 0$, then $I(\tau_0)>0$ for
some $\tau_0\in (0, h_0)$. From this there easily follow (\ref{eq4.47b})
and (\ref{eq4.47c}). $\Box$

\bigskip

\begin{nonsec}{Integral of energy and allocation of finite
forms }\label{4.48a}
\end{nonsec}
There is another application of the above estimates of energy integrals
connected with a generalization of the classical Denjoy--Carleman--Ahlfors
theorem about the number of different asymptotic tracts of an entire
function of a given order. In the present case this theorem can be
interpreted as a statement concerning the connection between the number
of finite forms in the class $\WT_2$ defined on the manifold $\M$ and
the rate of growth of their energy integrals.

Let $\M$ be an $n$-dimensional noncompact Riemannian manifold with or
without boundary. We fix a locally Lipschitz exhaustion function
$h:\M\to (0,h_0)$, $0<h_0\le\infty$, of the manifold $\M$.

We assume that there are $L\ge 1$ mutually disjoint domains
$\O_1,\O_2,\ldots,\O_L$  on $\M$  such that $\O_i\cap\partial\M=\emptyset$ 
if the boundary $\partial\M$ is nonempty. We also assume that on each domain
$\O_i$ is given a differential form $Z_i$ with continuous
coefficients and the properties:

\smallskip
$\deg Z_i=k-1,\quad dZ\not\equiv 0$,

\smallskip
$dZ_i=w\in\WT_2$ with structure constants $p,\nu_1,\nu_2$, independent
of $i=1,2,\ldots,L$,

\smallskip
$Z_i$ satisfies on $\partial\O_i$ the zero boundary condition
(\ref{eq4.4}).

\smallskip
We define a form $Z$ on $\M$ by setting $Z|_{\O_i}=Z_i$ and $Z=0$ on
$\M\setminus\cup_{i=1}^L\O_i$.

According to Theorem \ref{4.20} each of the domains $\O_i$ has a noncompact
closure. Then by Theorem \ref{4.46} the "narrower" the intersection of the
domains $\O_i$ with $h$-spheres $\Sigma_h(t)$ for $t\to\infty$, the
higher is the rate of growth of the form $Z$. Below we shall consider
the Denjoy--Carleman--Ahlfors theorem as a statement on the connection
between the number $L$ of mutually disjoint domains $\O_i$ on $\M_0$ and
the rate of growth of the energy of the form $Z$ (or of the form $Z$
itself) with respect to an exhaustion function $h(m)$ of the manifold $\M$.
We shall prove that such a formulation of the problem contains,
in particular, the classical Denjoy--Carleman--Ahlfors problem for
holomorphic functions of the complex plane. In the case of harmonic
functions of ${\bf R}^n$ see \cite{HK} for the history of the problem.

\smallskip

%\begin{subsec}{}\label{4.49}
%\end{subsec}
 We next introduce some necessary notation. We consider an open
subset $D\subset\M$ with a noncompact closure and we assume that the
restriction of the form $Z$ to $D$ satisfies condition (\ref{eq4.4}).

The function $h|_D:D\to(0,\infty)$ is an exhaustion function of $D$.
We fix an $h$-ball $B_h(\tau)$. Considering the variational problem
(\ref{eq4.34}) for the class of forms $Z$, satisfying the boundary condition
(\ref{eq4.4}) on $D_1$ we define the characteristic
$$
\epsilon(t;D)=\epsilon(t;\F_0)
$$
where $\F_0$ is defined in \ref{subsbas} and in (\ref{eq4.3}).

Following \cite{MIK1} we introduce the $N$-mean
\eqb
E(t;N)=\inf{1\over N}\sum\limits_{k=1}^N\epsilon(t;D_k)
\label{eq4.50}
\eqe
where the infimum is taken over all decompositions of $D$
into $N$ nonintersecting open sets $D_1,D_2,\ldots,D_N$ with noncompact closures.

We record the following simple result.
\smallskip

\begin{lem}{}\label{4.51}
 Let $D_1\subset D_2$ be arbitrary  open subsets of $\M$ 
with noncompact closures. Then
\eqb
\epsilon(t;D_2) \le\epsilon(t;D_1)\label{eq4.52}
\eqe
and
\eqb
\label{eq4.53}
\epsilon(t;\M)\le E(t;N)\,,\quad N\ge 1.
\eqe
\end{lem}

{\bf Proof.} It is enough to observe that each pair of forms $Z,Z_0$
admissible for the variational problem (\ref{eq4.34}) for the set $D_1$ is
also admissible for this problem for the set $D_2$.

 From (\ref{eq4.52}) we get (\ref{eq4.53}).
$\Box$

We next derive a more general assertion about the monotonicity of
$N$-means.

\begin{lem}{}\label{4.54}
For arbitrary $N>1$ we have
\eqb
E(t;N+1)\ge E(t;N). \label{eq4.55}
\eqe
\end{lem}

{\bf Proof.} We consider an arbitrary family of open subsets $\{D_k\},
k=1,2,\ldots,$ $N+1$, admissible for  the infimum in
(\ref{eq4.50}). It is not difficult to see that
$$
{1\over N+1}\sum_{k=1}^{N+1}\epsilon(t;D_k)={1\over N+1}\sum_{k=1}^{N+1}
\Bigl({1\over N}\sum_{j=1,j\neq k}^{N+1}\epsilon(t;D_j)\Bigr).
$$
Because
$$
{1\over N}\sum_{j=1,j\neq k}^{N+1}\epsilon(t;D_k)\ge E(t;N),
$$
we see that
$$
{1\over N+1}\sum_{k=1}^{N+1}\epsilon(t;D_k)\ge E(t;N)
$$
and the lemma is proved.
$\Box$
\bigskip

The next theorem provides a solution to the aforementioned problem
concerning the connection between the number $L$ of finite forms
on $\M$, and the rate of growth of the total energy of these forms
or the sum of their $L^p$-norms on an $h$-ball $B_h(\tau)$.

\begin{thm}{}\label{4.56}
Suppose that the manifold $\M$ satisfies the properties listed in the
beginning of this subsection and that for some $N=1,2,\ldots$ 
\eqb
\int\limits^{h_0}E(t,N) dt=\infty.\label{eq4.57}
\eqe
If the differential form $Z$, $dZ\in\WT_2(\M)$, is such that
\eqb
\label{eq4.58a}
\liminf_{\tau\to h_0}\int\limits_{h(m)<\tau}|dZ|^p\,*\1\,\exp\left\{-
\nu_1 \int\limits_{\tau_0}^{\tau}E(t;N)\,dt\right\}=0,
\eqe
or
\eqb
\liminf_{\tau',\tau''\to h_0\atop{0<\tau'<\tau''<h_0}}\;\;\;\int\limits_{\tau'<h(m)<\tau''}|\nabla h|\,|Z\wedge
\theta|*\1\exp\left\{-\nu_1\int\limits_{\tau_0}^{\tau'}
E(t;N)dt\right\}=0,
\label{eq4.58}
\eqe
or
\eqb
\liminf_{\tau',\tau''\to h_0\atop{0<\tau'<\tau''<h_0}}
\;\;\;\int\limits_{\tau'<h(m)<\tau''}|\nabla h|^p|Z|^p*
\1\exp\left\{-\nu_1\int\limits_{\tau_0}^{\tau'}
E(t;N)dt\right\}=0\,,
\label{eq4.59}
\eqe
then $L<N$.
\end{thm}

{\bf Proof.} We assume that there exists $N$ mutually nonintersecting
domains $\O_1,\O_2,\ldots,\O_N$ on the set $\M_0$ and the forms $Z_i$
defined on $\O_i$ with above properties. We denote by
$$
d(\O_k)=\inf_{m\in\O_k}h(m),\quad d=\max_{1\le k\le N} d(\O_k).
$$

Fix $\tau_0>d$. Using the inequality (\ref{eq4.37}) from Theorem \ref{4.35}
for an arbitrary $k=1,2,\ldots,N$ and a.e. $0<\tau_0<\tau'<h_0$ we have
$$
I_k(\tau_0)\exp\left\{\nu_1
\int\limits_{\tau_0}^{\tau'}\epsilon_k(t)\,dt\right\}\le I_k(\tau'),
$$
where
$$
I_k(\tau')=\int\limits_{\O_k\cap B_h(\tau')}|dZ_k|^p*\1,\quad
\epsilon_k(\tau')=\epsilon(\tau';\O_k).
$$
Adding these inequalities, we get
$$
\min_{1\le k\le N}I_k(\tau_0)\sum_{k=1}^N
\exp\left\{\nu_1
\int\limits_{\tau_0}^{\tau'}\epsilon_k(t)\,dt\right\}\le I(\tau'),
$$
where
$$
I(\tau')=\int\limits_{B_h(\tau')}|dZ|^p*\1.
$$
Applying the arithmetic--geometric mean inequality 
$$
{1\over N}\,\sum_{k=1}^N
\exp\left\{\nu_1
\int\limits_{\tau_0}^{\tau'}\epsilon_k(t)\,dt\right\} \ge
\prod_{k=1}^N\exp\left\{{\nu_1\over N}
\int\limits_{\tau_0}^{\tau'}\epsilon_k(t)\,dt\right\}
$$
we get
$$
\min_{1\le k\le N}I_k(\tau_0)
N\exp\left\{
{\nu_1\over N}
\int\limits_{\tau_0}^{\tau'}{1\over N}\sum_{k=1}^N\epsilon_k(t)\,dt\right\}\le I(\tau').
$$

The domains $\O_1,\O_2,\ldots,\O_N$ are nonintersecting. Therefore for
all $\tau_0<t<h_0$ we have
$$
{1\over N}\sum_{k=1}^N\epsilon_k(t)\ge E(t;N).
$$
The preceding inequality gives now
$$
\min_{1\le k\le N}I_k(\tau_0)N\exp\left\{\nu_1\int\limits_{\tau_0}^{\tau'} E(t;N)dt\right\}
\le  I(\tau').
$$

For the estimation of the integral $I(\tau')$ we use the inequalities
(\ref{eq4.41}) and (\ref{eq4.43b}) and obtain
$$
\min_{1\le k\le N}I_k(\tau_0)\le C_1\int\limits_{\tau'<h(m)<\tau''}
|\nabla h||Z\wedge\theta|*\1\exp\left\{-\nu_1\int\limits_{\tau_0}^{\tau'} E(t;N)dt\right\}
$$
or
$$
\min_{1\le k\le N} I_k(\tau_0)\le C_2\int\limits_{\tau'<h(m)<\tau''}
|\nabla h|^p|Z|^p*\1\;\exp\left\{-\nu_1
\int\limits_{\tau_0}^{\tau'} E(t;N)dt\right\}
$$
where 
$$
C_1={p\over{\tau''-\tau'}},\quad C_2=\left({{p\nu_2}\over{(\tau''-\tau')\nu_1}}\right).
$$

On the basis of the conditions (\ref{eq4.57})--(\ref{eq4.59}) imposed on
the form $Z$, for some $k$, $1\le k\le N$, we have $I_k(\tau_0)=0$. Thus
$dZ_k(m)\equiv 0$ on $\O_k\cap B_h(\tau_0)$. Because we have chosen
$\tau_0 >d$ arbitrarily, we can conclude that at least on one of the
components $\O_k$, $dZ_k\equiv 0$. Contradiction.
$\Box$

\bigskip

%\newpage

 \noindent
 {\bf Martio }\\
 Department of Mathematics and Statistics\\
 00014 University of Helsinki\\
 FINLAND\\
 Email: {\tt martio@cc.helsinki.fi}\\
 %Fax: +358-9-19151400\\

 \medskip

 \noindent
{\bf Miklyukov}\\
Mathematics Department\\
Volgograd State University\\
2 Prodolnaya 30\\
Volgograd 400062 \\
RUSSIA\\
E-mail:{\tt miklyuk@hotmail.com}\\
% Fax: +7-8442-471608

\noindent
{\bf Vuorinen}\\
Department of Mathematics\\
FIN-20014 University of Turku \\
FINLAND\\
E-mail: {\tt vuorinen@utu.fi}\\
%Fax: 358-2-333 6595

\end{document}